\documentclass[11 pt]{amsart}
\usepackage{amssymb}
\usepackage{amsmath}
\usepackage{amsfonts}
\usepackage{graphicx}
\usepackage{amsthm}
\usepackage{enumerate}
\usepackage[mathscr]{eucal}
\usepackage{verbatim}

 \theoremstyle{plain}
\newtheorem{theorem}{Theorem}[section]
\newtheorem{lemma}[theorem]{Lemma}

\numberwithin{equation}{section}
\theoremstyle{plain}

\theoremstyle{remark}

\def\bbR{{\mathbb {R}}}
\def\bbZ{{\mathbb {Z}}}

\def\cA{\mathcal A}

\def\cJ{\mathcal J}

\begin{document}

\date{October, 2008}

\title
{A uniform estimate for Fourier restriction to simple curves}

\author[]
{Daniel M. Oberlin}

\address
{D. M.  Oberlin \\
Department of Mathematics \\ Florida State University \\
 Tallahassee, FL 32306}
\email{oberlin@math.fsu.edu}

\thanks{The author was supported in part by NSF grant DMS-0552041.
}



\begin{abstract}
We prove a uniform Fourier extension-restriction estimate for a certain class of curves
in $\bbR ^d$.
\end{abstract}

\maketitle

\section{Introduction}

Let $\gamma$ be a curve in $\bbR^d$ given by 
\begin{equation}\label{ineq0}
\gamma (t) =\Big(t,\frac{t^2}{2} ,\dots ,\frac{t^{d-1}}{(d-1)!},\phi (t)\Big)
\end{equation}
where $\phi\in C^{(d+2)}(a,b)$ and $\phi ^{(j)}(t)> 0$ for $t\in (a,b)$, 
$j=0,1,2,\dots ,d+2$. Such curves are termed {\it simple} in \cite{DM1}.
Write $\omega(t)$ for $\phi^{(d)}(t)$. The purpose of this note is to prove a uniform Fourier
extension-restriction theorem for affine arclength measure on curves \eqref{ineq0}:

\begin{theorem}\label{theorem2} Suppose $\lambda$ is the measure on $\gamma$ 
given by $d\lambda =\omega (t)^{2/(d^2+d)}dt$. If $1\leq p<d+2$ and $\frac{1}{p}+\frac{d(d+1)}{2}\frac{1}{q}=1$, then there is $C(p,d)$ such that the following estimate holds: 
\begin{equation*}
\|\widehat{f\,d\lambda}\|_q \leq C(p,d)\,\|f\|_{L^p (\lambda)}.
\end{equation*}
\end{theorem}

\noindent If $d=2$ this is just the theorem in \cite{O1}, a result originally established
in slightly more generality, but with a more complicated proof, by Sj\"olin in \cite{S}. This note
is the result of an attempt to apply the method from \cite{O1} in higher dimensions. Theorem \ref{theorem2}, which
is an immediate consequence 
of Theorems \ref{theorem1} and \ref{theorem3} below,
 is somewhat analogous to the result of \cite{DW} 
on two fronts: it is a direct consequence of a geometric inequality (Theorem \ref{theorem1} 
here) combined with a fairly simple argument (Theorem \ref{theorem3} here), and 
its range of exponents $p$ is the (probably suboptimal) range obtained by Christ \cite{Ch}. 
For a better range of $p$ when $\omega$ satisfies a certain auxiliary condition, see 
Theorem 1.1 in \cite{BOS2}. For some of the
history of the problem of restricting Fourier transforms to curves, see \cite{BOS1}.

\begin{theorem}\label{theorem1} There exists $C=C(d)$ such that the estimate 
\begin{equation}\label{mainest}
\int_a^b \Big(\int_{\{t_1 \leq t_i <b \}}\chi_F \big(\gamma (t_1)+\cdots +\gamma (t_d )\big)
\Big[\prod_{i=2}^d \omega (t_i )\Big]^{2/(d^2 +d)}dt_2 \cdots dt_d \Big)^{(d+2)/2}
 \omega (t_1 )^{2/(d^2 +d)}\,dt_1
 \end{equation}
 \begin{equation*}
 \leq C(d) \,m_d (F)
\end{equation*}
holds for Borel $F\subset \bbR^d$.
\end{theorem}

\noindent The next result is Theorem 7 in \cite{O}.
\begin{theorem}\label{theorem3} Suppose that $\lambda$ is a nonnegative Borel
measure on $\bbR^d$ satisfying the inequality 
$$
\int\Bigl(\int_{\{\tau (y_1 )\leq \tau (y_i )\}}{\chi}_F
(y_1 +y_2 +\cdots +y_m )\, d\lambda (y_2 )
\cdots d\lambda (y_m )\Bigr)^{{m+2}\over{2}}d\lambda (y_1 )
\leq c\ m_d (F)$$
for some nonnegative integer $m\geq 3$, some real-valued Borel function
$\tau$ on supp$(\lambda )$, and all Borel $F\subseteq {\Bbb R}^d$. Then the adjoint
restriction estimate 
$$
\|{\widehat {fd\lambda}}\|_q \leq C(c,p) \ \|f\|_{L^p (\lambda)} $$
holds whenever ${1\over p}+{{m(m+1)}\over 2}{1 \over q} =1$ and $1\leq p
<m+2$.
\end{theorem}

\noindent The next section contains the proof of Theorem \ref{theorem1},
while \S3 contains the proofs of certain lemmas used in \S 2.

\noindent{\bf Acknowledgment:} The author would like to express  his gratitude to Jong-Guk Bak and Andreas Seeger:
this note is a byproduct of the collaboration that led to \cite{BOS1} and \cite{BOS2}.

\section{Proof of Theorem \ref{theorem1}}

The idea for the proof of Theorem \ref{theorem1}, from
\cite{O}, is to regard \eqref{mainest} as an $L^{(d+2)/2, 1}\rightarrow
L^{(d+2)/2}$ estimate for a certain operator $T$ and to establish \eqref{mainest} by establishing the dual 
$L^{(d+2)/d}\rightarrow L^{(d+2)/d,\infty}$ estimate for $T^\ast$.
Let $J(t_1 ,\dots ,t_d )$ be the absolute value of the Jacobian determinant for the map 
\begin{equation*}
(t_1 ,\dots t_d )\mapsto \gamma (t_1)+\cdots +\gamma (t_d ).
\end{equation*}
Then
\begin{equation*}
\langle f ,T^\ast g \rangle =
\langle Tf,g\rangle=
\end{equation*}
\begin{equation*}
\int_a^b \Big(\int_{\{t_1 \leq t_i <b \}}f \big(\gamma (t_1)+\cdots +\gamma (t_d )\big)
\Big[\prod_{i=2}^d \omega (t_i )\Big]^{2/(d^2 +d)}dt_2 \cdots dt_d \Big)g(t_1 )\,
 \omega (t_1 )^{2/(d^2 +d)}\,dt_1=
\end{equation*}
\begin{equation*}
\int_{\{a<t_1 \leq t_i <b \}}\Big(\frac{g(t_1 )}{J(t_1 ,\dots ,t_d)}
{\Big[\prod_{i=1}^d \omega (t_i )\Big]^{2/(d^2 +d)}}\Big)
f \big(\gamma (t_1)+\cdots +\gamma (t_d )\big){J(t_1 ,\dots ,t_d)}\,dt_d \cdots dt_1 \,
\end{equation*}
so that the desired $L^{(d+2)/d}\rightarrow L^{(d+2)/d,\infty}$ estimate for nonnegative $g$,
\begin{equation*}
\lambda ^{(d+2)/d} m_d \{T^\ast g \geq \lambda\}\leq C(d) \int_a^b g(t_1)^{(d+2)/d}
\omega (t_1 )^{2/(d^2 +d)}\, dt_1 \ ( \lambda >0) ,
\end{equation*} 
will follow from the estimate
\begin{equation}\label{ineq2}
\lambda ^{(d+2)/2}
\int
\chi_{\widetilde E} (t_1 ,\dots ,t_d )\, J(t_1 ,\dots ,t_d )\, dt_d \cdots dt_1
\leq 
C(d) \int_a^b g(t_1)^{(d+2)/d}\,
\omega (t_1 )^{2/(d^2 +d)}\, dt_1 ,
\end{equation}
where 
\begin{equation*}
\widetilde E =\Big\{(t_1 ,\dots ,t_d ): a<t_1 \leq\cdots\leq t_d <b ,\, \frac{g(t_1 )}{J(t_1 ,\dots ,t_d)}
{\Big[\prod_{i=1}^d \omega (t_i )\Big]^{2/(d^2 +d)}}\geq \lambda \Big\}.
\end{equation*}
(The change of variables implicit in this argument can be justified as in \cite{DM2}, p. 549.)
By absorbing $\lambda$ into $g$ we can assume $\lambda =1$.
Thus \eqref{ineq2} will follow from integrating the inequality
\begin{equation}\label{ineq4}
\int\chi_{\widetilde E} (t_1 ,\dots ,t_d )\, J(t_1 ,\dots ,t_d )\, dt_d \cdots dt_2
\leq 
C(d)\, g(t_1)^{(d+2)/d}\,
\omega (t_1 )^{2/(d^2 +d)}
\end{equation}
with respect to $t_1$. Lemma 2.3 in \cite{BOS2} shows that there is a nonnegative function $\psi =\psi (u;t_1 ,\dots ,t_d )$ supported in $[t_1 ,t_d ]$ such that
\begin{equation}\label{ineq4.5}
J(t_1 ,\dots ,t_d )=\int_{t_1}^{t_d}\omega (u)\,\psi (u;t_1 ,\dots ,t_d )\, du 
\end{equation}
and so \eqref{ineq4} will follow from the inequality
\begin{equation}\label{ineq5}
\int\chi_{\widetilde E} (t_1 ,\dots ,t_d )\, J(t_1 ,\dots ,t_d )\, dt_d \cdots dt_2
\leq 
C(d)\, c^{(d+2)/d}\,
\omega (t_1 )^{2/(d^2 +d)},
\end{equation}
to hold for $c>0$ and $t_1 \in (a,b)$, where now 
\begin{equation*}
\widetilde{E}=\Big\{(t_1 ,\dots ,t_d ): a<t_1 \leq\cdots\leq t_d <b ,\,\int_{t_1}^{t_d}\omega (u)\,
\psi (u;t_1 ,\dots ,t_d )\, du \leq c
{\Big[\prod_{i=1}^d \omega (t_i )\Big]^{2/(d^2 +d)}} \Big\}.
\end{equation*}
Homogeneity allows absorbing $c$ into $\omega$, so we can assume $c=1$. 
With $t_1 >a$ fixed, then, and with 
\begin{equation}\label{ineq7}
E=\Big\{(t_2 ,\dots ,t_d ): t_1\leq t_2 \leq\cdots\leq t_d <b ,\,J(t_1 ,\dots ,t_d )\leq
{\Big[\prod_{i=1}^d \omega (t_i )\Big]^{2/(d^2 +d)}}\Big \},
\end{equation}
inequality \eqref{ineq5}, and so \eqref{mainest}, will follow from 
\begin{equation}\label{ineq9}
\int\chi_{ E} (t_2 ,\dots ,t_d )\, 
\Big[\prod_{i=2}^d \omega (t_i )\Big]^{2/(d^2 +d)}
\, dt_d \cdots dt_2
\leq 
C(d).
\end{equation}

To begin the proof of \eqref{ineq9}, 
let $J\subset\bbZ$ be an interval of integers such that 
$\{2^j :j\in J\}$ is the set of dyadic values assumed by $\omega$ on $(a,b)$. For each $j\in J$, choose 
$a_j \in (a,b)$ such that $\omega (a_j )=2^j$. If $J$ has a least element, say $j_{\min}$, we let $a_{j_{\min} -1}=a$
and append $j_{\min} -1$ to $J$. If $J$ has a greatest element, we make a similar accommodation. Then,
writing $I_j =[a_j ,a_{j+1})\cap (a,b)$, we obtain a partition $\{ I_j \}_{j\in J}$ of $(a,b)$.

Now, with $t_1 \in (a,b)$ fixed, say $t_1 \in I_{j_1}$, with $E$ 
as in \eqref{ineq7}, and for integers $j_2 \leq j_3 \leq\cdots \leq j_d$ in 
$J$ with $j_2 \geq j_1$, we set
\begin{equation*}
E_{j_2 \cdots j_d}\dot=\{(t_2 ,\dots ,t_d)\in E:(t_2 ,\dots ,t_d )\in I_{j_2}\times\cdots
\times I_{j_d}\}.
\end{equation*}
The desired estimate \eqref{ineq9} will follow from 
\begin{equation}\label{ineq15}
\sum_{j_2  \geq j_1}\cdots \sum_{j_d \geq j_{d-1}} \big( 2^{j_2 +\cdots +j_d}\big)^{2/(d^2 +d)}
m_{d-1}\big(
E_{j_2 \cdots j_d }
\big)
\leq C(d).
\end{equation}
To establish \eqref{ineq15} it is enough to show that, for each $(j_2 ,\dots ,j_d )$ 
figuring in the sum in \eqref{ineq15},
we have
\begin{equation}\label{ineq16}
 \big( 2^{j_2 +\cdots +j_d}\big)^{1/(d-1)}
m_{d-1}\big(
E_{j_2 \cdots j_d}
\big)^{d/2}\leq C(d)\, \big( 2^{j_1 +j_2 +\cdots +j_d}\big)^{2/(d^2 +d)}.
\end{equation}
In fact, some algebra shows that \eqref{ineq16} is equivalent to
\begin{equation*}
 \big( 2^{j_2 +\cdots +j_d}\big)^{2/(d^2 +d)}
m_{d-1}\big(
E_{j_2 \cdots j_d}
\big)
\leq C(d)\, 2^{4j_1 /(d^3 +d^2 )}2^{-4(j_2 +\cdots +j_d)/[(d^3 +d^2 )(d-1)]}
\end{equation*}
and so, given \eqref{ineq16},  \eqref{ineq15} follows by summing a geometric series.

Moving towards the proof of \eqref{ineq16}, fix $(j_2 ,\dots ,j_d )$. 
In what follows we will often write $j(l)$ instead of $j_l$.
Let ${p_1}<{p_2}<\cdots {p_{k-1}}$ be the indices $i$
in $\{1,2 ,\dots ,d-1\}$ for which $j({i+1})-j(i) \geq 2$ and set $p_0 =0$ and $p_k =d$. 
Define $\ell_1 ,\dots ,\ell_k$ by $\ell_n =p_n -p_{n-1}-1$ ($n=1,\dots ,k$), and observe that $\ell_1 +\cdots +\ell_k =d-k$.
Then
\begin{equation*}
\{j(1), j(2) ,\dots ,j(d) \}=
\end{equation*}
\begin{equation*}
\{j({p_0 +1}),j(p_0 +2),\dots ,j({p_1});j({p_1 +1}),\dots ,j({p_2});\dots ;
j({p_{k-1}+1}),\dots ,j({p_k})\}
\end{equation*}
where if $j(i)$ and $j({i+1})$ are separated by a semicolon then $j({i+1})-j(i) \geq 2$ and 
otherwise $0\leq  j({i+1})-j(i) \leq 1$.
Next we construct $k$ subintervals $J_n$ of $(a,b)$ by setting, for $n=1,\dots ,k$, 
\begin{equation}\label{ineq17.5}
J_n =I_{j({p_{n-1}}+1)}\cup I_{j({p_{n-1}+2)}} \cup \cdots \cup I_{j({p_{n})}} 
\end{equation}
so that, recalling the definition of $I_j$, the endpoints of $J_n$ are $c_n \doteq a_{j({p_{n-1}+1)}}$ and $d_n 
\doteq a_{j({p_n}) +1}$.
Note that $c_1 <d_1 <c_2 <d_2 <\cdots  <c_k <d_k$ (see \eqref{ineq18.5} below) and that
if $(t_2 ,\dots ,t_d )\in E_{j_2 \dots  j_d}$ then
\begin{equation}\label{ineq18}
c_n \leq t_{p_{n-1}+1}\leq t_{p_{n-1}+2}\leq \cdots \leq t_{p_{n}}\leq d_n .
\end{equation}
We will need the facts that if $n=2,\dots k$, then
\begin{equation}\label{ineq18.5}
d_{n-1}\leq a_{j(p_{n-1}+1)-1}<c_n
\end{equation}
and
\begin{equation}\label{ineq19}
c_n - a_{j({p_{n-1}+1})-1}\gtrsim d_n -c_n .
\end{equation}
(Through this note, the constants implied by symbols like
$\lesssim$ can easily be checked to depend only on $d$.) To see \eqref{ineq19}, note that because $\omega (a_j )=2^j$ and $\omega '$ is nondecreasing we have
\begin{equation*}
(a_{j+1}-a_j )\,\omega '(a_j )\leq \int_{a_{j}}^{a_{j+1}}\omega '(u)\,du =2^j =2  \int_{a_{j-1}}^{a_{j}}
\omega '(u)\,du \leq 2 \,(a_{j}-a_{j-1} )\,\omega '(a_{j} )
\end{equation*}
so that $(a_{j+1}-a_j ) \leq 2 \,(a_{j}-a_{j-1} )$ and therefore 
\begin{equation}\label{ineq20}
(a_{j+p}-a_{j+p-1} ) \leq 2^p \,(a_{j}-a_{j-1} ).
\end{equation}
Now, by definition of $p_{n-1}$, $j(p_{n-1})+1\leq j(p_{n-1}+1)-1$, and so $a_{j(p_{n-1}+1)-1}$
lies between $d_{n-1}=a_{j(p_{n-1})+1}$ and $c_n =a_{j({p_{n-1}+1)}}$ in the sense of \eqref{ineq18.5}. 
Also, according to \eqref{ineq17.5},  $J_n =(c_n ,d_n )$ is (up to endpoints) the union of no more 
than $d$ intervals $I_j =[a_j ,a_{j+1})$. By choice of the $p_n$, each interval but the first in the union in \eqref{ineq17.5} 
is either identical to or contiguous to the one on its left. Since the first of these intervals is
$(a_{j({p_{n-1}+1)}},a_{j({p_{n-1}+1})+1})$,
\eqref{ineq20} implies \eqref{ineq19}.

We now outline the proof of \eqref{ineq16}, beginning with a lemma (the proofs of the lemmas will be given in \S\ref{proofs}): 

\begin{lemma}\label{lemma1}
Suppose $t_1 <\cdots <t_d$ and $(\alpha_i ,\beta_i )\subset (t_i ,t_{i+1})$. Write
$\Delta_i =\beta_i -\alpha_i$ and suppose $$f=\sum_{i=1}^{d-1}c_i \, \chi_{(\alpha_i ,\beta_i )}$$ where $c_i \geq 0$.
Fix $p\in\{1,\dots ,d-1\}$. Suppose $$\{e_i : i=1 ,2,\dots d-1, i\not= p\}=\{1,2,\dots ,d-2\}.$$ 
Then 
\begin{equation*}
\label{ineq10}
\int_{t_1}^{t_d} f(u)\,\psi (u;t_1 ,\dots ,t_d )\, du \gtrsim
c_p \,\Delta_p^{d-1}
\prod_{\genfrac{}{}{0pt}{}{1 \leq i \leq d-1}{i \neq p}} 
\Delta_i ^{e_i}.
\end{equation*}
\end{lemma}

\noindent With $t_1$ and $j_2 ,\dots ,j_d$ fixed and 
with $(t_2 ,\dots ,t_d )\in E_{j_2 \dots  j_d}$, we will apply Lemma \ref{lemma1}
to a collection $\mathcal I$ of  intervals $(\alpha_i ,\beta_i )$ specified as follows: for $n=1,\dots ,k$ and $i=p_{n -1}+1 ,\dots ,p_n -1$, the  $\ell_n$ intervals $(t_i ,t_{i+1})$ will be in $\mathcal I$; additionally, for $n=2 ,\dots ,k$, the intervals 
\begin{equation}
\label{ineq30.25}
\tilde{J}_n \dot= (a_{j(p_{n-1}+1)-1},c_n)
\subset (t_{p_{n-1}},t_{p_{n-1}+1})
\end{equation}
 will be in $\mathcal I$ (we set $\widetilde{J}_1 =\emptyset$). 
Observe that there are integers $m_1 \leq m_2 \leq \cdots \leq m_k$ such that 
\begin{equation}\label{ineq30.5}
\omega \sim 2^{m_n}\  \text{on}\ 
\tilde{J}_n \cup\Big(\cup_{i=p_{n-1}+1}^{p_n -1}(t_i ,t_{i+1})\Big) \dot=\cJ_n \ (n=1,\dots ,k).
\end{equation}
(This is true because, according to \eqref{ineq18},  $\cJ_n$ is 
contained in the union of at most $\ell_n +1\leq d$ contiguous intervals $I_j$,  
and $\omega\sim2^j$ on $I_j$.)
Then \eqref{ineq16} can be written 
\begin{equation}\label{ineq31}
2^{[\ell_1 m_1 +(\ell_2 +1)m_2+\cdots +(\ell_{k}+1)m_k ]/(d-1)}
m_{d-1}\big( E_{j_2 \cdots j_d} \big)^{d/2}
\lesssim 
\end{equation}
\begin{equation*}
2^{2[(\ell_1 +1)m_1 +(\ell_2 +1)m_2+\cdots +(\ell_{k}+1)m_k ]/(d^2+d)}.
\end{equation*}
Similarly, the inequality 
\begin{equation*}
J(t_1 ,\dots ,t_d )\leq\Big[\prod_{i=1}^d \omega (t_i )\Big]^{2/(d^2 +d)}
\end{equation*}
in \eqref{ineq7}
implies
\begin{equation}\label{ineq32}
J(t_1 ,\dots ,t_d )\lesssim 2^{2[(\ell_1 +1)m_1 +\cdots +(\ell_{k}+1)m_k ]/(d^2+d)}.
\end{equation}
Now, as we will see below, Lemma \ref{lemma1}, \eqref{ineq30.5}, \eqref{ineq32}, and 
$$
J(t_1 ,\dots ,t_d )=\int_{t_1}^{t_d}\omega (u)\, \psi (u;t_1 ,\dots ,t_d )\, du
$$
will yield certain estimates of the form
\begin{equation}\label{ineq34}
2^{m_n}\Delta_p^{d-1}
\prod_{\genfrac{}{}{0pt}{}{1 \leq i \leq d-1}{i \neq p}} 
\Delta_i^{e_i}\lesssim 
2^{2[(\ell_1 +1)m_1 +(\ell_2 +1)m_2+\cdots +(\ell_{k}+1)m_k ]/(d^2+d)}
\end{equation}
when $(\alpha_p ,\beta_p )\subset\cJ_n$.
A weighted geometric mean of these estimates \eqref{ineq34} will give
\begin{equation}\label{ineq35}
2^{[\ell_1 m_1 +(\ell_2 +1)m_2+\cdots +(\ell_{k}+1)m_k ]/(d-1)}\,
W_1^{d/(\ell_1 +1)}
\prod_{n=2}^k \Big( \rho_n^{(d+\ell_n )/2}W_n^{(d-1)/(\ell_n +1)}\Big)\lesssim
\end{equation}
\begin{equation*}
2^{2[(\ell_1 +1)m_1 +(\ell_2 +1)m_2+\cdots +(\ell_{k}+1)m_k ]/(d^2+d)},
\end{equation*}
where $\rho_n $ is the length of $\tilde{J}_n$
and with the $W_n$'s given by
\begin{equation*}
W_n =W (t_{p_{n-1}+1},\dots ,t_{p_n})
\end{equation*}
where, for $s_1 \leq\cdots \leq s_m$,
\begin{equation*}
W(s_1 ,\dots ,s_m )=\sup\big\{\prod_{i=1}^{m-1}(s_{i+1}-s_i )^{e_i}:
\{e_1 ,\dots ,e_{m-1}\}=\{1 ,\dots ,m-1\}\big\}.
\end{equation*}
Lemma \ref{lemma2} below will allow the choice of $(t_2 ,\dots ,t_d )\in E_{j_2 \dots j_d} $ such that
\begin{equation}\label{ineq36}
m_{d-1}\big( E_{j_2 \cdots j_d} \big)^{d/2}\lesssim 
W_1^{d/(\ell_1 +1)}
\prod_{n=2}^k \Big( \rho_n^{(d+\ell_n )/2}W_n^{(d-1)/(\ell_n +1)}\Big).
\end{equation}
With \eqref{ineq35} this will yield \eqref{ineq31}.

To give the details missing from the argument in the preceding paragraph we will  need a lemma
whose statement requires the introduction of some more notation:
$\cA_{d-1}$ will stand for the convex hull in $\bbR^{d-1}$ of the set of all permutations of the $({d-1})$-tuple
$(1, 2,\dots ,{d-1})$. Recall that $\ell_1 +\cdots +\ell_k+k-1={d-1}$. If $\ell_1 >0$,
$\cA _{d-1} '$ is defined to be  
the collection  of all permutations of $({d-1})$-tuples
\begin{equation*}
{(d-1)}\Big( \frac{1}{\ell_1}\big(1 ,\dots ,\ell_1 \big); \frac{1}{\ell_2 +1}\big(1 ,\dots ,\ell_2 \big);\dots ;
 \frac{1}{\ell_k +1}\big(1 ,\dots ,\ell_k \big); \underbrace{\frac{1}{2},\cdots ,\frac{1}{2}}_{k-1\;\rm times}\Big).
\end{equation*}
Note that if $k=1$, then $\cA _{d-1} ' =\cA _{d-1} $.
For $k\geq 2$, define 
$\cA_{d-1}  ''$ to be the collection  of all permutations of ${(d-1)}$-tuples
\begin{equation*}
{(d-1)}\Big( \frac{1}{\ell_1 +1}\big(1 ,\dots ,\ell_1 \big); \dots ;
 \frac{1}{\ell_k +1}\big(1 ,\dots ,\ell_k \big);\underbrace{\frac{1}{2},\cdots ,\frac{1}{2}}_{k-2\;\rm times};1\Big).
\end{equation*}
(To simplify the notation, and since no confusion will result from doing so,
we surpress the dependence of $\cA_{d-1} '$ and $\cA_{d-1} ''$ on $k$ and 
the $\ell_n$'s.) 

\begin{lemma}\label{lemma4} The inclusions $\cA_{d-1} ' ,\,\cA_{d-1} '' \subset \cA_{d-1}$ hold.
\end{lemma}

Moving towards \eqref{ineq35}, fix $n' \in\{2,\dots ,k\}$. We will show that
\begin{equation}\label{ineq41}
2^{m_{n'}}\prod_{n=1}^k W_n^{(d-1)/(\ell_n +1)}
\Big(
\prod_{\genfrac{}{}{0pt}{}{1 \leq n \leq k}{n \neq 1,n'}} 
\rho_n^{(d-1)/2}\Big)\rho_{n'}^{d-1} \lesssim
\end{equation}
\begin{equation*}
2^{2[(\ell_1 +1)m_1 +\cdots +(\ell_{k}+1)m_k ]/(d^2+d)}.
\end{equation*}
Recall the definitions of the intervals $(\alpha_i ,\beta_i )$, whose lengths $\Delta_i$ are the numbers 
$\rho_n$ ($n=2,\dots ,k$) along with the numbers 
$t_{i+1}-t_i$ for $i=p_{n-1}+1,\dots ,p_n -1$ and $n=1,\dots ,k$. Since
\begin{equation*}
W_n =\prod_{i=p_{n-1}+1}^{p_n -1}(t_{i+1}-t_i )^{e_i^n},
\end{equation*}
for some choice of $\{e_i^n  \}$ with 
$$
\{e_i^n \}_{i=p_{n-1}+1}^{p_n -1}=\{1,\dots ,\ell_n \} ,
$$
it follows that 
\begin{equation}\label{ineq41.5}
\prod_{n=1}^k W_n^{(d-1)/(\ell_n +1)}\Big(
\prod_{\genfrac{}{}{0pt}{}{1 \leq n \leq k}{n \neq 1,n'}} 
 \rho_n^{(d-1)/2}\Big)\rho_{n'}^{d-1}=
\prod_{i=1}^{d-1} \Delta_i ^{\sigma (i)},
\end{equation} 
where the vector $\sigma =\big(\sigma (i)\big)$ is in $\cA_{d-1} ''$ and where, if
$i_0$ is the index for which $(\alpha_{i_0},\beta_{i_0})=\widetilde{J}_{n'}$, then
$\sigma (i_0 )=d-1$ and so
$\Delta_{i_0}^{\sigma (i_0 )}=\rho_{n'}^{d-1}$.
By Lemma \ref{lemma4}, the vector $\sigma $ is a convex combination 
\begin{equation}\label{ineq42}
\sigma =\sum_q \lambda_q \tau_q
\end{equation}
of vectors $\tau_q =\big(\tau_q (i)\big)$, each of which is a permutation of $(1 ,\dots ,d-1)$.
Further, since $\sigma (i_0 )=d-1$, we have $\tau_q (i_0 )=d-1$ for each $q$. Now it follows from Lemma \ref{lemma1} that if the $c_i$'s are nonnegative, then
\begin{equation}\label{ineq42.5}
c_{i_0} \,\Delta_{i_0}^{d-1}
\prod_{\genfrac{}{}{0pt}{}{1 \leq i \leq d-1}{i \neq i_0}} 
\Delta_i^{\tau_q (i)}\lesssim
\int_{t_1}^{t_d} \sum_{i=1}^{d-1}c_i \, \chi_{(\alpha_i ,\beta_i)}(u)\,\psi (u;t_1 ,\dots ,t_d )\, du 
\end{equation}
for each $q$. From \eqref{ineq41.5},  \eqref{ineq42}, and \eqref{ineq42.5} it then follows that
\begin{equation*}
c_{i_0}\prod_{n=1}^k W_n^{(d-1)/(\ell_n +1)}\Big(
\prod_{\genfrac{}{}{0pt}{}{1 \leq n \leq k}{n \neq 1,n'}} 
\rho_n^{(d-1)/2}\Big)\rho_{n'}^{d-1}\lesssim
\int_{t_1}^{t_d} \sum_{i=1}^{d-1}c_i \, \chi_{(\alpha_i ,\beta_i)}(u)\,\psi (u;t_1 ,\dots ,t_d )\, du .
\end{equation*}
Since \eqref{ineq30.5} implies that $\omega\sim 2^{m_{n'}}$ on $\widetilde{J}_{n'}=
(\alpha_{i_0},\beta_{i_0})$, we have
\begin{equation}\label{ineq44}
2^{m_{n'}}\prod_{n=1}^k W_n^{(d-1)/(\ell_n +1)}
\Big(\prod_{\genfrac{}{}{0pt}{}{1 \leq n \leq k}{n \neq 1,n'}} 
 \rho_n^{(d-1)/2}\Big)\rho_{n'}^{d-1}\lesssim
\int_{t_1}^{t_d}2^{m_{n'}}\chi_{\widetilde{J}_{n'}}(u) \,\psi (u;t_1 ,\dots ,t_d )\, du \lesssim
\end{equation}
\begin{equation*}
\int_{t_1}^{t_d}\omega (u) \,\psi (u;t_1 ,\dots ,t_d )\, du =J(t_1 ,\dots ,t_d )\lesssim
 2^{2[(\ell_1 +1)m_1 +\cdots +(\ell_{k}+1)m_k ]/(d^2+d)}
\end{equation*}
by \eqref{ineq4.5} and \eqref{ineq32}. This is \eqref{ineq41}. Analogous to \eqref{ineq41} we will also need, in the case $\ell_1 >0$,
the estimate
\begin{equation}\label{ineq45}
2^{m_1}W_1^{(d-1)/\ell_1}
\prod_{n=2}^k W_n^{(d-1)/(\ell_n +1)}
\prod_{n=2}^k
\rho_n^{(d-1)/2}\ \lesssim
\end{equation}
\begin{equation*}
2^{2[(\ell_1 +1)m_1 +\cdots +(\ell_{k}+1)m_k ]/(d^2+d)}.
\end{equation*}
As before, 
\begin{equation*}
W_1^{(d-1)/\ell_1}
\prod_{n=2}^k W_n^{(d-1)/(\ell_n +1)}
\prod_{n=2}^k
\rho_n^{(d-1)/2}\ =
\prod_{i=1}^{d-1} \Delta_i ^{\sigma (i)}
\end{equation*} 
where now $\sigma$ is in $\cA_{d-1} '$. With $\sigma =\sum_q \lambda_q \tau_q$ 
as in \eqref{ineq42}, Lemma \ref{lemma1} gives
$$
2^{m_1} \,\prod_{i=1}^{d-1}\Delta_i^{\tau_q (i)}\lesssim
\int_{t_1}^{t_d} \sum_1^{d-1}2^{m_1} \, \chi_{(\alpha_i ,\beta_i)}(u)\,\psi (u;t_1 ,\dots ,t_d )\, du 
$$
for each $q$. This leads, as before, to
$$
2^{m_1}W_1^{(d-1)/\ell_1}
\prod_{n=2}^k W_n^{(d-1)/(\ell_n +1)}
\prod_{n=2}^k
\rho_n^{(d-1)/2}\ \lesssim
\int_{t_1}^{t_d} \sum_1^{d-1}2^{m_1} \, \chi_{(\alpha_i ,\beta_i)}(u)\,\psi (u;t_1 ,\dots ,t_d )\, du .
$$
Since $\omega\gtrsim 2^{m_1}$ on $[t_1 ,t_d ]$, \eqref{ineq45} follows as in \eqref{ineq44}.

Now \eqref{ineq35} will follow by considering a particular weighted geometric mean of the estimates 
\eqref{ineq45} and \eqref{ineq41}. In fact, given the computations
\begin{equation*}
\frac{\ell_1}{d-1}+\frac{\ell_2 +1 }{d-1}+\frac{\ell_3 +1 }{d-1}+\cdots +\frac{\ell_k +1 }{d-1} =1,
\end{equation*} 
\begin{equation*}
\frac{d-1}{\ell_1}\frac{\ell_1}{d-1}+\frac{d-1}{\ell_1 +1}\Big(\frac{\ell_2 +1 }{d-1}+\frac{\ell_3 +1 }{d-1}+\cdots +\frac{\ell_k +1 }{d-1}
\Big)=\frac{d}{\ell_1 +1},
\end{equation*}
and
\begin{equation*}
(d-1)\frac{\ell_{n'} +1}{d-1}+\frac{d-1}{2}\Big( 1-\frac{\ell_{n'} +1}{d-1}\Big)=\frac{d+\ell_{n'}}{2},\ {n'}=2,\dots ,k,
\end{equation*}
\eqref{ineq35} is an immediate consequence of \eqref{ineq41} and \eqref{ineq45}.

Now the proof of \eqref{ineq31} will be complete when we have explained how to choose $(t_2 ,\dots ,t_d )\in E_{j_2 \dots j_d}$ so that 
\eqref{ineq36} holds. We will need another lemma and 
some more notation: recall that
\begin{equation*}
c_1 <d_1 <c_2 <d_2 <\cdots <c_k <d_k .
\end{equation*}
Let $\delta _n =d_n -c_n$. 
Recall that $1\leq p_1 <p_2 <\cdots <p_k =d$, $p_0 =0$, $\ell_n =p_n -p_{n-1}-1$, and 
that 
\begin{equation}\label{ineq30}
c_n \leq t_{p_{n-1}+1}<\cdots <t_{p_n}\leq d_n
\end{equation}
for $n=1,\dots k$. 
With $t_{p_{n-1}+1}\in [c_n ,d_n )$, 
write $\bf{t}_n$ for an $\ell_n$-tuple 
$(t_{p_{n-1}+2},\dots ,t_{p_n})$ satisfying 
\eqref{ineq30},
and $\bf{t}$ for the $(\ell_1 +\cdots +\ell_k =d-k)$-tuple $(\bf{t}_1 ,\dots ,\bf{t}_k )$. 

\begin{lemma}\label{lemma2}
The inequality 
\begin{equation*}
m_{d-k} \big(
\{
{\bf{t}}:
 W_1^{d/(\ell_1 +1)}
\prod_{n=2}^k  W_n^{(d-1)/(\ell_n +1)}\leq \mu
 \}
\big)\lesssim
\mu^{2/d}\prod_{n=2}^k \delta_n^{\ell_n /d}
\end{equation*}
holds for $\mu >0$.
\end{lemma}

Writing 
$$
(t_2 ,\dots ,t_d )=(t_{p_1 +1},\dots ,t_{p_{k-1}+1};\bf{t})
$$
with {\bf t} as above, we use \eqref{ineq30} to choose 
$$
(t_{p_1 +1}',\dots ,t_{p_{k-1}+1}')\in \prod_{n=2}^k \ [c_n ,d_n ]
$$
such that
\begin{equation*}
m_{d-k} \big(
\{{\bf{t}}: (t_{p_1 +1}',\dots ,t_{p_{k-1}+1}';{\bf{t}})\in E_{j_2 \cdots j_d}\}
\big)
\geq \frac{m_{d-1}(E_{j_2 \cdots j_d})}{\prod_{n=2}^k \delta_n}.
\end{equation*}
Let $c_1 (d)>0$ be sufficiently small. Then if $\mu >0$ is such that
$$
\mu^{2/d}\prod_{n=2}^k \delta_n^{\ell_n /d} =c_1 (d)\frac{m_{d-1}(E_{j_2 \cdots j_d})}{\prod_{n=2}^k \delta_n},
$$
it follows from Lemma \ref{lemma2} that there is 
$$
(t_2 ,\dots ,t_d)\in E_{j_2 \cdots j_d}
$$ 
such that 
$$
 W_1^{d/(\ell_1 +1)}
\prod_{n=2}^k  W_n^{(d-1)/(\ell_n +1)}> \mu =
c_2 (d)\frac{m_{d-1}( E_{j_2 \cdots j_d})^{d/2}}{\prod_{n=2}^k \delta_n^{(d+\ell_n )/2}}.
$$
Recalling that $\rho_n$ is the length of $\tilde{J}_n$, so that 
$\delta_n =d_n -c_n \lesssim \rho_n$ by \eqref{ineq30.25} and \eqref{ineq19}, \eqref{ineq36} follows.

\section{Proofs of lemmas }\label{proofs}

\noindent{\bf{Proof of Lemma \ref{lemma1}.}}
The proof is by induction on $d$ and, since $\psi (u;t_1 ,t_2 )=\chi_{[t_1 ,t_2 ]}(u)$, the case 
$d=2$ is clear. 
Fix $p\in\{1,\dots ,d-1\}$ and then 
a choice of $\{e_i \}_{i\not= p}$ so that $\{e_i \}_{i\not= p}=\{1,\dots ,d-2\}$. 
Let ${q}$ satisfy $e_{{q}}=1$.
We will give the argument in the case ${q}<p$, the case ${q}>p$ 
being similar. Let $m_i$ be the midpoint of $[\alpha_i ,\beta_i ]$.
Define intervals $I_1 ,\dots I_{d-1}$ and $J_1 ,\dots J_{d-1}$ as follows:
\begin{equation*}
I_j =(\alpha_i ,m_i ),\, J_i =(m_i ,\beta_i ) \ \text{if}\, i<{q},
\end{equation*}
\begin{equation*}
I_{{q}} =(\alpha_{{q}} ,m_{{q}} ),\, J_{{q}} =\emptyset,\ \ \ \ \ \ \ \ \ \ \ \ \ \ \ \ \ \ 
\end{equation*}
\begin{equation*}
I_i =(m_i ,\beta_i ),\, J_i =(\alpha_i ,m_i ) \ \text{if}\, i>{q}.
\end{equation*}
We use the identity
\begin{equation*}
\psi (u;t_1 ,\dots ,t_d )=\int_{t_1}^{t_2}\cdots \int_{t_{d-1}}^{t_d}\psi (u;s_1 ,\dots ,s_{d-1} )\,
ds_1 \cdots ds_{d-1},
\end{equation*}
a consequence of the proof of Lemma 2.3 in \cite{BOS2},
and the fact that $\psi (u,s_1 ,\dots ,s_{d-1} )$ is supported on $[s_1 ,s_{d-1}]$ to write
\begin{equation}\label{ineq11}
\int_{t_1}^{t_d} f(u)\,\psi (u;t_1 ,\dots ,t_d )\, du
=\int_{t_1}^{t_2}\cdots \int_{t_{d-1}}^{t_d}
\int_{s_1}^{s_{d-1}}f(u)\,
\psi (u;s_1 ,\dots ,s_{d-1} )\,du\,
ds_1 \cdots ds_{d-1}.
\end{equation}
Then 
\begin{equation*}
\eqref{ineq11}\geq
\int_{I_1}\cdots \int_{I_{d-1}}
\int_{s_1}^{s_{d-1}}
\sum_{{\genfrac{}{}{0pt}{}{1 \leq i \leq d-1}{i \neq q}}}c_i \,\chi_{J_i}(u)
\,
\psi (u;s_1 ,\dots ,s_{d-1} )\,du\,
ds_1 \cdots ds_{d-1}.
\end{equation*}
Now if $s_i \in I_i$ for $i=1,\dots ,d-1$ then $J_i \subset (s_i ,s_{i+1})$ if $i<{q}$ and 
$J_i \subset (s_{i-1} ,s_{i})$ if $i>{q}$.
Thus, assuming the lemma for $d-1$ and noting that $\{e_i  -1\}_{i\not= p,{q}} =\{1,\dots ,d-3\}$, it follows that 
\begin{equation*}
\int_{s_1}^{s_{d-1}}
\sum_{{\genfrac{}{}{0pt}{}{1 \leq i \leq d-1}{i \neq q}}}c_i \,\chi_{J_i}(u)
\,
\psi (u;s_1 ,\dots ,s_{d-1} )\,du
\gtrsim 
c_p \,\Delta_p^{d-2}
\prod_{\genfrac{}{}{0pt}{}{1 \leq i \leq d-1}{i \neq p,q}}
\Delta_i^{e_i -1}
\end{equation*}
and so
\begin{equation*} 
\eqref{ineq11}\gtrsim c_p \,\big(\prod_{i=1}^{d-1}\Delta_i \big)\Delta_p^{d-2}
\prod_{\genfrac{}{}{0pt}{}{1 \leq i \leq d-1}{i \neq p,q}}
\Delta_i^{e_i -1}
=
c_p \,\Delta_p^{d-1}
\prod_{\genfrac{}{}{0pt}{}{1 \leq i \leq d-1}{i \neq p}}
\Delta_i^{e_i },
\end{equation*}
completing the proof of Lemma \ref{lemma1}.
\medskip

\noindent{\bf{Proof of Lemma \ref{lemma4}.}}
Lemma \ref{lemma4} is the statement that if 
$$\ell_1 +\cdots +\ell_k+k-1=r,$$ then {\bf (a)} if $\ell_1 >0$, $k\geq 2$,  and
$\cA _r '$ is the collection  of all permutations of $r$-tuples
\begin{equation}\label{ineq64.2}
r\Big( \frac{1}{\ell_1}\big(1 ,\dots ,\ell_1 \big); \frac{1}{\ell_2 +1}\big(1 ,\dots ,\ell_2 \big);\dots ;
 \frac{1}{\ell_k +1}\big(1 ,\dots ,\ell_k \big); \underbrace{\frac{1}{2},\cdots ,\frac{1}{2}}_{k-1\;\rm times}\Big),
\end{equation}
we have $\cA _r ' \subset\cA_r$ and {\bf (b)} if $k\geq 2$ and
$\cA_r  ''$ is the collection  of all permutations of $r$-tuples
\begin{equation*}
r\Big( \frac{1}{\ell_1 +1}\big(1 ,\dots ,\ell_1 \big); \dots ;
 \frac{1}{\ell_k +1}\big(1 ,\dots ,\ell_k \big);\underbrace{\frac{1}{2},\cdots ,\frac{1}{2}}_{k-2\;\rm times};1\Big),
\end{equation*}
we have $\cA_r '' \subset \cA_r$. We will show these inclusions by
induction on $k$. 
We require the following two facts, which we establish at the end of the proof of this lemma:
\begin{equation}\label{ineq64.9}
\Big(\frac{1}{\ell_{k-1}+1}\big(1,\dots ,\ell_{k-1}\big);
\frac{1}{\ell_{k}+1}\big(1,\dots ,\ell_{k}\big)\Big)\in\Big(\frac{1}{\ell_{k-1}+\ell_k +1} \Big) \cA_{\ell_{k-1}+\ell_k}
\end{equation}
and
\begin{equation}\label{ineq65}
\Big(\frac{1}{\ell_{k-1}+1}\big(1,\dots ,\ell_{k-1}\big);
\frac{1}{\ell_{k}+1}\big(1,\dots ,\ell_{k}\big);\frac{1}{2}\Big)\in\Big(\frac{1}{\ell_{k-1}+\ell_k +2}\Big) \cA_{\ell_{k-1}+\ell_k +1}.
\end{equation}

If $k=2$ and $\ell_1 >1$, a vector \eqref{ineq64.2} can be written
\begin{equation*}
r\Big( \frac{1}{\ell_1 }\big(1 ,\dots ,\ell_1 -1\big);1; \frac{1}{\ell_{2} +1}\big(1 ,\dots ,\ell_{2}\big);\frac{1}{2}
\Big),
\end{equation*}
and therefore, by \eqref{ineq65}, as a linear combination of permutations of vectors
\begin{equation*}
r\Big( \frac{1}{\ell_1 +\ell_2 +1}\big(1 ,\dots ,\ell_1 +\ell_2 \big);1\Big)\in\cA_r .
\end{equation*}
And if $k=2$ and $\ell_1 =1$, then \eqref{ineq64.2} can be written
\begin{equation*}
r\Big( 1; \frac{1}{\ell_{2} +1}\big(1 ,\dots ,\ell_{2}\big);\frac{1}{1+1}\big(1 \big)
\Big)
\end{equation*}
and therefore, by \eqref{ineq64.9}, as a linear combination of permutations of the vector
\begin{equation*}
r\Big( 1; \frac{1}{\ell_{2} +2}\big(1 ,\dots ,\ell_{2}+1\big)
\Big)\in\cA_r .
\end{equation*}
Thus {\bf (a)} holds for $k=2$.
The fact that {\bf (b)} holds for $k=2$,
follows similarly 
from \eqref{ineq64.9}. So assume that $k\geq 3$ and that {\bf (a)} and {\bf (b)} hold with $k-1$ in place of $k$.

To show that {\bf (b)} holds for $k$, fix a vector 
\begin{equation}\label{ineq65.5}
r\Big( \frac{1}{\ell_1 +1}\big(1 ,\dots ,\ell_1 \big);\dots ; \frac{1}{\ell_{k-1} +1}\big(1 ,\dots ,\ell_{k-1} \big);
 \frac{1}{\ell_k +1}\big(1 ,\dots ,\ell_k \big); \underbrace{\frac{1}{2},\cdots ,\frac{1}{2}}_{k-2\;\rm times};1\Big)
\end{equation}
in $\cA_r ''$. It follows from \eqref{ineq65} that
\eqref{ineq65.5} can be written as a convex combination of permutations of vectors 
\begin{equation}\label{ineq66}
r\Big( \frac{1}{\ell_1 +1}\big(1 ,\dots ,\ell_1 \big);\dots ; \frac{1}{\ell_{k-2} +1}\big(1 ,\dots ,\ell_{k-2} \big);
 \frac{1}{\ell_{k-1}+\ell_k +2}\big(1 ,\dots ,\ell_{k-1}+\ell_k +1\big);
  \underbrace{\frac{1}{2},\cdots ,\frac{1}{2}}_{k-3\;\rm times};1\Big),
\end{equation}
and our induction assumption implies that each permutation of \eqref{ineq66} is in $\cA_r$. This establishes 
{\bf (b)} for $k$.

To see that {\bf (a)} holds for $k$, note that the argument above shows that a vector 
\begin{equation*}
r\Big( \frac{1}{\ell_1 }\big(1 ,\dots ,\ell_1 \big);\dots ; \frac{1}{\ell_{k-1} +1}\big(1 ,\dots ,\ell_{k-1} \big);
 \frac{1}{\ell_k +1}\big(1 ,\dots ,\ell_k \big); \underbrace{\frac{1}{2},\cdots ,\frac{1}{2}}_{k-1\;\rm times}\Big)
\end{equation*}
of the form \eqref{ineq64.2} can be written as a convex combination of permutations of vectors 
\begin{equation*}
r\Big( \frac{1}{\ell_1}\big(1 ,\dots ,\ell_1 \big);\dots ; \frac{1}{\ell_{k-2} +1}\big(1 ,\dots ,\ell_{k-2} \big);
 \frac{1}{\ell_{k-1}+\ell_k +2}\big(1 ,\dots ,\ell_{k-1}+\ell_k +1\big);
  \underbrace{\frac{1}{2},\cdots ,\frac{1}{2}}_{k-2\;\rm times}\Big),
\end{equation*}
and so, by the induction assumption, is in $\cA_r$. 

It remains to establish \eqref{ineq64.9} and \eqref{ineq65}.
We require an alternate description of $\cA_r$. Let $\cA_r^*$ be the set 
\begin{equation*}
\Big\{ (a_1 ,\dots ,a_r ):\sum_{j=1}^r a_{j}= \frac{r(r+1)}{2}\ \text{and}\, 
\sum_{j \in E} a_j \geq\frac{|E|(|E|+1)}{2} \ \text{if}\, E\subset\{1,\dots ,r\}\Big\}.
\end{equation*}
We want to show
that $\cA_r =\cA_r^*$, and  it is enough to show that each extreme point of the convex set
$\cA_r^*$ is a permutation of $(1 ,\dots ,r)$. 
So assume that $(a_1 ,\dots ,a_{r}) $ is an extreme point of
$\cA_r^*$. Without loss of generality we may also assume $a_1 \leq a_2 \leq\cdots\leq a_r$.

Our first step will be to show that 
\begin{equation}\label{ineq63.1}
a_1 <a_2 <\cdots <a_d .
\end{equation}
To this end, assume that $a_s =a_{s+1}$ for some $s\in\{1,\dots ,d-1\}$. We will show that if
either $s\in E,s+1\notin E$ or  $s+1\in E,s\notin E$, then  
\begin{equation}\label{ineq63.2}
\sum_{j\in E}a_j >\frac{|E|(|E|+1)}{2}.
\end{equation}
If \eqref{ineq63.2} holds, it will follow that there is $\delta>0$ such
the vector obtained from 
$(a_1 ,\dots ,a_{r}) $ by replacing $a_s$ and $a_{s+1}$ by 
$a_s +\eta$ and $a_{s+1}-\eta$ is in $\cA_r$ whenever $|\eta |<\delta$.
This implies that $(a_1 ,\dots ,a_{r}) $ is not extreme.
To show \eqref{ineq63.2} we begin with an observation:
\begin{equation}\label{ineq63.3}
\text{if}\ a_k =a_{k+1},\ \text{then}\ \sum_{j=1}^k a_j >\frac{k(k+1)}{2}. 
\end{equation}
(To see \eqref{ineq63.3}, observe that the assumption 
$\sum_{j=1}^k a_j =\frac{k(k+1)}{2}$
and the inequality $\sum_{j=1}^{k+1} a_j \geq\frac{(k+1)(k+2)}{2}$ together imply that 
$a_{k+1}\geq k+1$ and so $a_k \geq k+1$ if $a_k =a_{k+1}$. Then 
$$
\sum_{j=1}^k a_j =a_k +\sum_{j=1}^{k-1} a_j \geq (k+1)+\frac{(k-1)k}{2},
$$
contradicting $\sum_{j=1}^k a_j =\frac{k(k+1)}{2}$.) Returning to 
\eqref{ineq63.2}, we will write $k=|E|$ and consider three cases:

\noindent{\bf Case I} ($k<s$):  Here we will show that if $s\in E$ or $s+1\in E$ then 
\eqref{ineq63.2} holds.
Assuming $s\in E$, it follows that
$$
\frac{k(k+1)}{2}\leq \sum_{j=1}^k a_j =\sum_{j=1}^{k-1} a_j +a_k \leq
\sum_{j=1}^{k-1} a_j +a_s \leq \sum_{j\in E}a_j .
$$ 
Thus if \eqref{ineq63.2} fails, then $a_k =a_{k+1}$ (since $a_k =a_s$) and
$ \sum_{j=1}^k a_j =k(k+1)/2$, contradicting \eqref{ineq63.3}.
The case $s+1\in E$ is similar.

\noindent{\bf Case II} ($k=s$): We will just observe that if $a_s =a_{s+1}$,
then \eqref{ineq63.2} holds. In fact, since
$$
\sum_{j\in E}a_j \geq \sum_{j=1}^k a_j ,
$$
this follows immediately from  \eqref{ineq63.3}.

\noindent{\bf Case III} ($k>s$): We will show that if 
either $s\notin E$ or $s+1\notin E$ then \eqref{ineq63.2} holds.
Write $E=\{j_1 ,j_2 ,\dots ,j_k \}$ with $j_1 <\cdots <j_k$. 
Since $k>s$ and either $s\notin E$ or $s+1\notin E$, it follows that 
$j_k >k$. Then, since $1\leq j_1 ,\dots ,k\leq j_k$, the inequality
$$
\frac{k(k+1)}{2}\leq \sum_{l=1}^k a_l 
 \leq \sum_{l=1}^k a_{j_l} =\sum_{j\in E}a_j
$$
shows that if \eqref{ineq63.2} fails then 
$a_k =a_{j_k}$ and $\sum_{j=1}^k a_j =k(k+1)/2$, again resulting in a contradiction of 
\eqref{ineq63.3}. Thus \eqref{ineq63.2}, and so \eqref{ineq63.1}, are established.

Now suppose that $(a_1 ,\dots ,a_{r}) $ is extreme and \eqref{ineq63.1} holds.
If $a_j \geq j $ for $j=1 ,\dots ,r$ then the condition 
$$
\sum_{j=1}^r a_j  =r(r+1)/2
$$
forces $(a_1 ,\dots ,a_{r})=(1 ,\dots ,r)$.  
But if we have $a_j <j$ for any $j$, we can choose $t$ such that
$a_t <t$ and that
$a_j \geq j $ for $j=1 ,\dots ,t-1$ (the condition $a_1\geq 1$ implies that $t>1$).
Since 
\begin{equation}\label{ineq64}
\sum_{j=1}^t a_j  \geq t(t+1)/2
\end{equation}
we can choose $s$ with $s\leq t-1$, $a_s >s$, and $a_j=j$ if $s<j<t$. 
Thus 
$$
(a_1 ,\dots ,a_r )=(a_1 ,\dots ,a_{s-1},a_s  ,{s+1},\dots ,{t-1},a_t  ,a_{t+1},\dots ,a_r ).
$$
It follows from \eqref{ineq64} and $a_t <t$ that 
$$
\sum_{j=1}^p a_j  > p(p+1)/2
$$
for $p=s,\cdots ,t-1$. Thus there is $\delta >0$ such that if $|\eta |<\delta$, then 
\begin{equation}\label{ineq64.1}
(a_1 ,\dots ,a_{s-1},a_s +\eta ,{s+1},\dots ,{t-1},a_t -\eta ,a_{t+1},\dots ,a_r )\in \cA_r^* ,
\end{equation}
where we have used the fact that \eqref{ineq63.1} implies that the
entries of the vector in \eqref{ineq64.1} are nondecreasing if $\delta$ is small enough.
Then $(a_1 ,\dots ,a_r )$ cannot be an extreme point of $\cA_r^*$. Thus $a_j \geq j$
for $j=1 ,\dots ,r$ and so $(a_1 ,\dots ,a_r )=(1,\dots ,r)$ as desired.

We return to the proofs for \eqref{ineq64.9} and \eqref{ineq65}. Since 
$\cA_{\ell_{k-1}+\ell_k}=\cA_{\ell_{k-1}+\ell_k}^*$,
\eqref{ineq64.9} will follow from checking that if $m\leq \ell_{k-1}$ and
$n\leq\ell_k$, then 
\begin{equation}\label{ineq70}
\frac{\ell_{k-1}+\ell_k +1}{\ell_{k-1}+1}\frac{m(m+1)}{2}+\frac{\ell_{k-1}+\ell_k +1}{\ell_{k}+1}\frac{n(n+1)}{2}\geq
\frac{(n+m)(n+m+1)}{2}.
\end{equation}
This inequality is equivalent to the inequality
\begin{equation}\label{ineq72}
m(m+1)\ell_{k}(\ell_k +1)+n(n+1)\ell_{k-1}(\ell_{k-1}+1)\geq 2mn(\ell_{k-1}+1)(\ell_{k}+1).
\end{equation}
And the easily checked inequality
\begin{equation*}
\sqrt{m(m+1)n(n+1)\ell_{k-1} (\ell_{k-1} +1)\ell_k (\ell_k +1)}\geq mn(\ell_{k-1}+1)(\ell_{k}+1)
\end{equation*}
shows that \eqref{ineq72} follows from the inequality between arithmetic and geometric means.

Similarly, to show that \eqref{ineq65} holds it is enough to check that if $m\leq \ell_{k-1}$ and
$n\leq\ell_k$, then the inequalities 
\begin{equation}\label{ineq73}
\frac{\ell_{k-1}+\ell_k +2}{\ell_{k-1}+1}\ \frac{m(m+1)}{2}+\frac{\ell_{k-1}+\ell_k +2}{\ell_{k}+1}\ \frac{n(n+1)}{2}\geq
\frac{(n+m)(n+m+1)}{2}
\end{equation}
and
\begin{equation}\label{ineq74}
\frac{\ell_{k-1}+\ell_k +2}{\ell_{k-1}+1}\ \frac{m(m+1)}{2}+\frac{\ell_{k-1}+\ell_k +2}{\ell_{k}+1}\ 
\frac{n(n+1)}{2}+\frac{\ell_{k-1}+\ell_k +2}{2}\geq
\end{equation}
\begin{equation*}\
\frac{(n+m+1)(n+m+2)}{2}
\end{equation*}
hold. Since \eqref{ineq73} follows from \eqref{ineq70}, it is enough to establish \eqref{ineq74}.
Now inequality \eqref{ineq70} is equivalent to 
\begin{equation*}
\frac{m(m+1)}{2(\ell_{k-1}+1)}+\frac{n(n+1)}{2({\ell_{k}+1})}\geq
\frac{(n+m)(n+m+1)}{2({\ell_{k-1}+\ell_k +1})},
\end{equation*}
while \eqref{ineq74} is equivalent to 
\begin{equation*}
\frac{m(m+1)}{2(\ell_{k-1}+1)}+\frac{n(n+1)}{2({\ell_{k}+1})}+\frac{1}{2}\geq
\frac{(n+m+1)(n+m+2)}{2({\ell_{k-1}+\ell_k +2})}=
\end{equation*}
\begin{equation*}
\Big(\frac{(n+m)(n+m+1)}{2}+(n+m+1)\Big)
\Big(\frac{1}{\ell_{k-1}+\ell_k +1}-\frac{1}{(\ell_{k-1}+\ell_k +1)(\ell_{k-1}+\ell_k +2)}\Big).
\end{equation*}
Thus it is enough to show that 
\begin{equation*}
\frac{1}{2}+\frac{(n+m)(n+m+1)}{2(\ell_{k-1}+\ell_k +1)(\ell_{k-1}+\ell_k +2)}
+\frac{n+m+1}{(\ell_{k-1}+\ell_k +1)(\ell_{k-1}+\ell_k +2)}\geq \frac{n+m+1}{\ell_{k-1}+\ell_k +1}.
\end{equation*}
This is equivalent to the inequality 
\begin{equation*}
(\ell_{k-1}+\ell_k +1)+\frac{(n+m+1)(n+m+2)}{\ell_{k-1}+\ell_k +2}\geq
2(n+m+1)
\end{equation*}
which follows from the arithmetic-geometric mean inequality and the easily checked 
\begin{equation*}
(n+m+1)(\ell_{k-1}+\ell_k +2)\leq (\ell_{k-1}+\ell_k +1)(n+m+2),
\end{equation*}
itself a consequence of $m\leq \ell_{k-1}$, $n\leq \ell_k$.

\noindent{\bf{Proof of Lemma \ref{lemma2}.}}
Recall that $c_n \leq t_{p_{n-1}+1}<\cdots <t_{p_n}\leq d_n$
and that $W_n =W (t_{p_{n-1}+1},\dots ,t_{p_n})$. It follows that $W_n \lesssim \delta_n^{\ell_n (\ell_n +1)/2}$
and so  
\begin{equation*}
m_{d-k} \big(\big\{{\bf{t}}:W_1^{d/(\ell_1 +1)}\prod_{n=2}^k  W_n^{(d-1)/(\ell_n +1)}\leq \mu\big\}\big)\lesssim
\end{equation*}
\begin{equation}\label{ineq50}
\sum_{2^{p_2}\leq  \delta_2^{\ell_2 (\ell_2 +1)/2}}\cdots \sum_{2^{p_k}\leq  \delta_k^{\ell_k (\ell_k +1)/2}}
\prod_{n=2}^k \Big|
\big\{ W_n \leq 2^{p_n}\big\}\Big|\cdot
\Big|\big\{ W_1 \leq\big(\mu /\prod_{n=2}^k 2^{p_n (d-1) /(\ell_n +1)}
\big)^{(\ell_1 +1)/d}\big\}\Big|,
\end{equation}
where 
$$
|\{W_n \leq \lambda\} |=m_{\ell_n}\big(\{(t_{p_{n-1}+2},\dots ,t_{p_n}):
W(t_{p_{n-1}+1},t_{p_{n-1}+2},\dots ,t_{p_n})\leq \lambda
\}\big).
$$
We will need the estimate      
\begin{equation}\label{ineq60}
|\{W_n \leq \lambda\} |\lesssim \lambda^{2/(\ell_n +1)}.
\end{equation}
To show \eqref{ineq60} and, more generally, to show that
\begin{equation}\label{ineq61}
m_p \big(\{(s_1 ,\dots ,s_p ):0\leq s_1 \leq\cdots\leq s_p :W(0,s_1 ,\dots ,s_p ) \leq \lambda\} \big)\leq C(p)\,\lambda^{2/(p +1)},
\end{equation}
we will argue by induction on $p$. 
(Inequality \eqref{ineq61} is an analog of a result, Proposition 2.4 (i) in \cite{BOS1}, from \cite{DM1} and \cite{DM2}.) 
The case $p=1$ is clear, so assume that \eqref{ineq61} 
holds. If we make the change of variable 
$$
u_1 =s_1 ,u_2 =s_2 -s_1 ,\dots ,u_p =s_p -s_{p-1},
$$
then 
$$
m_p \big(\{(s_1 ,\dots ,s_p ):0\leq s_1 \leq\cdots\leq s_p ,\, W(0,s_1 ,\dots ,s_p ) \leq \lambda\} \big)= 
$$
$$
p!\,m_p \big(\{(u_1 ,\dots ,u_p ):0\leq u_1 \leq\cdots\leq u_p ,\, \prod_{j=1}^p (u_j )^j \leq \lambda \}\big).
$$
Now
\begin{equation*}
m_{p+1} \big( \{(u_1 ,\dots ,u_{p+1} ):0\leq u_1 \leq\cdots\leq u_{p+1} ,\, \prod_{j=1}^{p+1} (u_j )^j \leq \lambda \}\big)=
\end{equation*}
\begin{equation*}
\int_0^{\infty}m_p \big(\{(u_1 ,\dots ,u_p ):0\leq u_1 \leq\cdots\leq u_p 
\leq u_{p+1} ,\, \prod_{j=1}^p u_j^j \leq {\lambda}/{ (u_{p+1})^{p+1}} \}
 \big)\,du_{p+1}.
\end{equation*}
Thus the estimate
$$
m_p \big(\{(u_1 ,\dots ,u_p ):0\leq u_1 \leq\cdots\leq u_{p+1} ,\, \prod_{j=1}^p (u_j )^j 
\leq {\lambda}/{ (u_{p+1})^{p+1}} \}\big)
\leq
$$
$$
C(p)\, \min\{(u_{p+1})^{p}, 
 ({\lambda}/{ (u_{p+1})^{p+1}})^{2/(p+1}\},
$$
a consequence of \eqref{ineq61}, shows that \eqref{ineq61} holds with $p$ replaced by $p+1$.

Now, using \eqref{ineq60} and some algebra, we have
\begin{equation*}
\eqref{ineq50}\lesssim \mu^{2/d}
\sum_{2^{p_2}\leq  \delta_2^{\ell_2 (\ell_2 +1)/2}}\cdots \sum_{2^{p_k}\leq  \delta_k^{\ell_k (\ell_k +1)/2}}
\prod_{n=2}^k  2^{2p_n /[d(\ell_n +1 ])}\lesssim \mu^{2/d}\prod_{n=2}^k \delta_n^{\ell_n /d}.
\end{equation*}
This gives the desired conclusion and therefore completes the proof of Lemma \ref{lemma2}.

\end{document}